\documentclass[12pt]{article}
\usepackage{amssymb}
\usepackage{latexsym}

\newtheorem{definition}{Definition}[section]
\newtheorem{theorem}[definition]{Theorem}
\newtheorem{lemma}[definition]{Lemma}
\newtheorem{corollary}[definition]{Corollary}
\newtheorem{remark}[definition]{Remark}

\newtheorem{problem}[definition]{Problem}
\newtheorem{note}[definition]{Note}

\newtheorem{proposition}[definition]{Proposition}

\def\K{\mathbb F}

\def\Z{\mathbb Z}

\begin{document}

\title{\bf 
$q$-Inverting
pairs of linear transformations\\
and the $q$-tetrahedron
algebra
}
\author{
Tatsuro Ito{\footnote{
Department of Computational Science,
Faculty of Science,
Kanazawa University,
Kakuma-machi,
Kanazawa 920-1192, Japan
}}
$\;$ and
Paul Terwilliger{\footnote{
Department of Mathematics, University of
Wisconsin, 480 Lincoln Drive, Madison WI 53706-1388 USA}
}}
\date{}

\maketitle
\begin{abstract}
As part of our  
study of the $q$-tetrahedron algebra $\boxtimes_q$
we introduce the notion of a {\it $q$-inverting pair}.
Roughly speaking, this is a pair of
invertible semisimple linear
transformations on a finite-dimensional vector space,
each of which acts on the eigenspaces of the other
according to a certain rule.
Our main result is
a bijection 
between the following two sets:
(i) the isomorphism classes of finite-dimensional irreducible 
$\boxtimes_q$-modules of type 1;
(ii) the isomorphism classes of $q$-inverting
pairs.
\medskip

 \end{abstract}

\section{Introduction}

\noindent Throughout this paper $\K$ denotes an algebraically
closed field.
We fix a nonzero $q\in \K$  that is not a root of 1.

\medskip
\noindent 
The $q$-tetrahedron algebra $\boxtimes_q$ 
was introduced
in \cite{qtet} as part of the ongoing investigation
of the Leonard pairs
\cite{hartwig},
\cite{nom2},
\cite{nom2a},
\cite{Rosengren2},
\cite{LS99},
  \cite{qSerre}, 
   \cite{LS24},
   \cite{conform},
    \cite{lsint},
\cite{qrac},
\cite{aw},
\cite{vidunas}
and  tridiagonal pairs
\cite{hasan},
\cite{hasan2},
\cite{TD00},
\cite{shape},
\cite{tdanduq},
\cite{NN}, 
\cite{nomura},
\cite{nom3},
\cite{nom4},
\cite{pasc}.
The algebra $\boxtimes_q$ is a unital associative
$\K$-algebra, and infinite-dimensional
as a vector space over $\K$. We defined 
$\boxtimes_q$ 
by generators and relations.
As explained in \cite{qtet},
   $\boxtimes_q$ can be viewed 
as a $q$-analog of the three-point loop
algebra $\mathfrak{sl}_2 \otimes_{\K} \K\lbrack t,t^{-1},(t-1)^{-1}\rbrack$
($t$ indeterminate).
The algebra $\boxtimes_q$ is related to 
the quantum group
$U_q(\mathfrak{sl}_2)$
\cite[Proposition 7.4]{qtet},
the 
$U_q(\mathfrak{sl}_2)$ loop algebra
\cite[Proposition 8.3]{qtet},
and the positive part of
$U_q({\widehat{\mathfrak{sl}}}_2)$ \cite[Proposition 9.3]{qtet}.
In \cite{qtet} we described the finite-dimensional 
irreducible 
$\boxtimes_q$-modules.
From this description and from \cite[Section 2]{NN}
there emerges 
a characterization of the finite-dimensional irreducible
$\boxtimes_q$-modules
in terms of 
a certain kind of tridiagonal pair
said to be {\it $q$-geometric}. For notational
convenience, in the present paper we will refer to
this as a {\it $q$-tridiagonal pair}.
As we will review in Section 3,
the following two sets are in bijection:
(i) the isomorphism classes of
finite-dimensional irreducible $\boxtimes_q$-modules
of type $1$;
(ii) the isomorphism classes of $q$-tridiagonal pairs.

\medskip
\noindent 
In the present paper we give
a second characterization of the 
finite-dimensional irreducible $\boxtimes_q$-modules,
this time using
a linear algebraic object called 
a {\it $q$-inverting pair}.
Roughly speaking, 
this is a pair of  invertible  semisimple  
linear transformations on a finite-dimensional
vector space, each of which acts on
the eigenspaces of the other according to a certain rule
that we find attractive. 
Our main result is 
a bijection between the following two sets:
(i) the isomorphism classes of finite-dimensional irreducible
$\boxtimes_q$-modules of type 1; 
(ii) the isomorphism classes of $q$-inverting
pairs.

\medskip
\noindent The plan for the paper is as follows.
In Section 2
we recall the algebra $\boxtimes_q$ and discuss its
finite-dimensional irreducible modules.
In Section 3  
we review the notion of a $q$-tridiagonal 
pair, and show how these objects are related to the 
finite-dimensional irreducible $\boxtimes_q$-modules.
In Section 4 we introduce the notion of a $q$-inverting
pair, and 
discuss how these objects  
are related to the finite-dimensional irreducible $\boxtimes_q$-modules.
 Theorem \ref{thm:easy} and
Theorem \ref{thm:hard} are the main results of the paper;
Sections 5--8 are devoted to their proofs.
In Section 9 we give some suggestions for further
research.

\section{The algebra $\boxtimes_q$}

\noindent
In this section we recall the
$q$-tetrahedron algebra and discuss
its finite-dimensional
irreducible modules.
We will use the following notation. Let $\Z_4=\Z/4\Z$
denote the cyclic group of order 4. Define
\begin{eqnarray*}
\lbrack n \rbrack_q = \frac{q^n-q^{-n}}{q-q^{-1}}
\qquad \qquad n=0,1,2\ldots
\end{eqnarray*}

\begin{definition}
\label{def:qtet}
\rm \cite[Definition 6.1]{qtet}
Let $\boxtimes_q$ denote the unital associative $\K$-algebra that has
generators 
\begin{eqnarray*}
\lbrace x_{rs}\;|\; r,s \in \Z_4,\;s-r=1 \;\mbox{or} \;s-r=2\rbrace
\end{eqnarray*}
and the following relations:
\begin{enumerate}
\item[{\rm (i)}]  
For $r,s\in \Z_4$ such that $s-r=2$,
\begin{eqnarray}
x_{rs}x_{sr} = 1.
\label{eq:t1}
\end{eqnarray}
\item[{\rm (ii)}]  
For $r,s,t\in \Z_4$ such that the pair $(s-r,t-s)$ is one of
$(1,1), (1,2), (2,1)$,
\begin{eqnarray}
\frac{qx_{rs}x_{st}-q^{-1}x_{st}x_{rs}}{q-q^{-1}}=1.
\label{eq:t2}
\end{eqnarray}
\item[{\rm (iii)}]  
For $r,s,t,u\in \Z_4$ such that $s-r=t-s=u-t=1$,
\begin{eqnarray}
\label{eq:qserre}
x_{rs}^3x_{tu} -
\lbrack 3 \rbrack_q
x_{rs}^2x_{tu}x_{rs} +
\lbrack 3 \rbrack_q
x_{rs}x_{tu}x_{rs}^2- 
x_{tu}x_{rs}^3=0.
\end{eqnarray}
\end{enumerate}
We call $\boxtimes_q$ the {\it $q$-tetrahedron algebra}.
\end{definition}
\begin{note}\rm
The equations (\ref{eq:qserre}) are the cubic $q$-Serre
relations 
\cite{charp}.
\end{note}
\noindent
We make some observations.

\begin{lemma}
\label{lem:rho}
There exists an $\K$-algebra automorphism $\rho$ of $\boxtimes_q$
that sends each generator $x_{rs}$ to $x_{r+1,s+1}$.
Moreover 
 $\rho^4=1$.
\end{lemma}

\begin{lemma}
\label{lem:flip}
There exists an $\K$-algebra automorphism of 
 $\boxtimes_q$
that sends each generator $x_{rs}$ to $-x_{rs}$.
\end{lemma}

\noindent We comment on the 
$\boxtimes_q$-modules.
Let $V$ denote a finite-dimensional irreducible 
$\boxtimes_q$-module. 
By 
\cite[Theorem 12.3]{qtet}
 there exist an
integer $d\geq 0$ and a scalar $\varepsilon \in \lbrace 1,-1\rbrace$
such that for each generator $x_{rs}$ the action on
$V$ is semisimple with eigenvalues $\lbrace \varepsilon
q^{d-2i}\;|\;0\leq i\leq d\rbrace$.
We call $d$ the {\it diameter} of $V$.
We call $\varepsilon $ the {\it type} of $V$.
Replacing each generator $x_{rs}$ by $\varepsilon x_{rs}$ the
type becomes 1.

%
%
%

\section{$q$-Tridiagonal pairs}

\medskip
\noindent In this section we recall the
notion of a {\it $q$-tridiagonal pair}
and discuss how these objects are related to
the finite-dimensional irreducible $\boxtimes_q$-modules.

\medskip
\noindent We will use the following notation.
Let $V$ denote a vector space over 
$\K$ with
finite positive dimension.
Let $\lbrace s_i\rbrace_{i=0}^d$
denote a finite sequence consisting of positive
integers whose sum is the dimension of $V$.
By a {\it decomposition of $V$
of shape 
$\lbrace s_i\rbrace_{i=0}^d$
} 
we mean
a sequence $\lbrace V_i\rbrace_{i=0}^d$
of subspaces of
$V$ such that 
$V_i$ has dimension $s_i$ for
$0 \leq i \leq d$ and 
\begin{eqnarray*}
V=\sum_{i=0}^d V_i \qquad \qquad (\mbox{direct sum}).
\end{eqnarray*}
We call $d$ the {\it diameter} of the decomposition.
For $0 \leq i \leq d$ we call $V_i$ the {\it $i^{th}$ component}
of the decomposition.
For notational convenience we define
$V_{-1}=0$ and $V_{d+1}=0$.
By the {\it inversion} of the decomposition
$\lbrace V_i\rbrace_{i=0}^d$ 
we mean the decomposition
$\lbrace V_{d-i}\rbrace_{i=0}^d $.

\begin{definition}
\rm
(\cite[Definition 1.1]{TD00}, \cite[Definition 2.6]{NN})
\label{def:tdq}
Let $V$ denote a vector space over $\K$ with finite positive
dimension. By a {\it $q$-tridiagonal pair on $V$}
we mean an ordered pair of linear transformations
$A:V\to V$ and $A^*:V \to V$ that satisfy (i)--(iii) below.
\begin{enumerate}
\item[{\rm (i)}]  
There exists a decomposition 
$\lbrace V_i\rbrace_{i=0}^d$ 
of $V$ such that
\begin{eqnarray*}
(A-q^{d-2i}I)V_i = 0  \qquad \qquad (0 \leq i \leq d),
\label{eq:td1}
\\
A^* V_i \subseteq V_{i-1} + V_i + V_{i+1}
\qquad \qquad (0 \leq i \leq d).
\label{eq:td2}
\end{eqnarray*}
\item[{\rm (ii)}]  
There exists a decomposition 
$\lbrace V^*_i\rbrace_{i=0}^{\delta}$ 
of
$V$ such that
\begin{eqnarray*}
(A^*-q^{\delta-2i}I)V^*_i = 0  \qquad \qquad (0 \leq i \leq \delta),
\label{eq:tds1}
\\
A V^*_i \subseteq V^*_{i-1} + V^*_i + V^*_{i+1}
\qquad \qquad (0 \leq i \leq \delta).
\label{eq:tds2}
\end{eqnarray*}
\item[{\rm (iii)}] There does not exist a subspace $W\subseteq V$
such that $AW \subseteq W$ and $A^*W \subseteq W$, other 
than $W=0$ and $W=V$.
\end{enumerate}
We say the pair $A,A^*$
is {\it over} $\K$. We call $V$
the {\it underlying vector space}.
\end{definition}

\begin{note} \rm
According to a common notational convention
$A^*$ denotes the conjugate transpose of $A$. We are
not using this convention. In a $q$-tridiagonal pair
$A,A^*$ 
the linear transformations $A$ and
$A^*$ are arbitrary subject
to 
(i)--(iii) above.
\end{note}

\begin{note}
\rm
The integers $d$ and $\delta$ from
Definition
\ref{def:tdq} are equal \cite[Lemma 4.5]{TD00};
we call this common value the {\it diameter} of
 the pair.
\end{note}

\noindent We now recall the notion of {\it isomorphism}
for $q$-tridiagonal pairs.

\begin{definition} \rm
Let $A,A^*$ and 
$A',A^{*\prime}$ denote $q$-tridiagonal pairs over $\K$.
By an {\it isomorphism of $q$-tridiagonal pairs} from 
$A,A^*$ to
$A',A^{*\prime}$
we mean a vector space isomorphism 
$\sigma$ from the vector space
underlying
$A,A^*$ to
the vector space underlying
$A',A^{*\prime}$
such that 
$\sigma A = A'\sigma$ and
$\sigma A^* = A^{*\prime} \sigma$.
We say that $A,A^*$ and
$A', A^{*\prime}$ are {\it isomorphic}
whenever there exists an isomorphism of $q$-tridiagonal
pairs from
 $A,A^*$ to
$A',A^{*\prime}$.
\end{definition}

\noindent Our results concerning
  $q$-tridiagonal pairs
and $\boxtimes_q$-modules 
are contained in the following two theorems and subsequent
remark.

\begin{theorem}
(\cite[Theorem 2.7]{NN},
\cite[Theorem 10.3]{qtet})
\label{thm:tdeasy}
Let $V$ denote a finite-dimensional irreducible
$\boxtimes_q$-module of type 1. Then the generators
$x_{01}, x_{23}$ act on $V$ as a $q$-tridiagonal pair.
\end{theorem}

\begin{theorem}
(\cite[Theorem 2.7]{NN}),
\cite[Theorem 10.4]{qtet})
\label{thm:tdhard}
Let $V$ denote a vector space over $\K$ with
finite positive dimension and let $A,A^*$ denote
a $q$-tridiagonal pair on $V$.
Then there exists a unique $\boxtimes_q$-module
structure on $V$ such that $x_{01}, x_{23}$ act
on $V$ as $A,A^*$ respectively.
This module structure is irreducible and type 1.
\end{theorem}

\begin{remark} 
\rm 
Combining Theorem \ref{thm:tdeasy} and
Theorem
\ref{thm:tdhard} we get a bijection between the
following two sets:
\begin{enumerate}
\item[{\rm (i)}]  
the isomorphism classes of finite-dimensional irreducible
$\boxtimes_q$-modules of type 1;
\item[{\rm (ii)}]  
the isomorphism classes of $q$-tridiagonal pairs.
\end{enumerate}
\end{remark}

\section{$q$-Inverting pairs}
\medskip
\noindent 
In this section we introduce the notion of
a {\it $q$-inverting pair} and
discuss how these objects are related to the finite-dimensional
irreducible $\boxtimes_q$-modules. This section
contains our main results.

\begin{definition}
\label{def:inv}
\rm
Let $V$ denote a vector space over $\K$ with finite positive
dimension. By a $q$-{\it inverting pair on $V$}
we mean an ordered pair of invertible linear transformations
$K:V\to V$ and $K^*:V \to V$ that satisfy (i)--(iii) below.
\begin{enumerate}
\item[{\rm (i)}]  
There exists a decomposition 
$\lbrace V_i\rbrace_{i=0}^d$ 
of $V$ such that
\begin{eqnarray}
(K-q^{d-2i}I)V_i = 0  \qquad \qquad (0 \leq i \leq d),
\label{eq:k1}
\\
K^* V_i \subseteq V_0 + V_1 + \cdots + V_{i+1}
\qquad \qquad (0 \leq i \leq d),
\label{eq:k2}
\\
(K^*)^{-1} V_i \subseteq V_{i-1} + V_i + \cdots + V_d
\qquad \qquad (0 \leq i \leq d).
\label{eq:k3}
\end{eqnarray}
\item[{\rm (ii)}]  
There exists a decomposition 
$\lbrace V^*_i\rbrace_{i=0}^{\delta}$ 
of
$V$ such that
\begin{eqnarray}
(K^*-q^{\delta-2i}I)V^*_i = 0  \qquad \qquad (0 \leq i \leq \delta),
\label{eq:ks1}
\\
K V^*_i \subseteq V^*_{i-1} + V^*_i + \cdots + V^*_\delta
\qquad \qquad (0 \leq i \leq \delta),
\label{eq:ks2}
\\K^{-1} V^*_i \subseteq V^*_0 + V^*_1 + \cdots + V^*_{i+1}
\qquad \qquad (0 \leq i \leq \delta).
\label{eq:ks3}
\end{eqnarray}
\item[{\rm (iii)}] There does not exist a subspace $W\subseteq V$
such that $KW \subseteq W$ and $K^*W \subseteq W$, other 
than $W=0$ and $W=V$.
\end{enumerate}
We say the pair $K,K^*$
is {\it over} $\K$. We call $V$
the {\it underlying vector space}.
\end{definition}

\begin{note} \rm
According to a common notational convention
$K^*$ denotes the conjugate transpose of $K$. We are
not using this convention. In a $q$-inverting pair
$K,K^*$ the linear transformations $K$ and
$K^*$ are arbitrary subject
to (i)--(iii) above.
\end{note}

\begin{note} \rm
The integers $d$ and $\delta$ from
Definition
\ref{def:inv} turn out to be equal; we will show this in
Lemma \ref{lem:dvsdel}.
\end{note}

\noindent We now define the notion of {\it isomorphism}
for $q$-inverting pairs.

\begin{definition} \rm
Let $K,K^*$ and 
$K',K^{*\prime}$ denote $q$-inverting pairs over $\K$.
By an {\it isomorphism of $q$-inverting pairs} from 
$K,K^*$ to
$K',K^{*\prime}$
we mean a vector space isomorphism 
$\sigma$ from the vector space
underlying
$K,K^*$ to
the vector space underlying
$K',K^{*\prime}$
such that 
$\sigma K = K'\sigma$ and
$\sigma K^* = K^{*\prime} \sigma$.
We say $K,K^*$ and
$K', K^{*\prime}$ are {\it isomorphic}
whenever there exists an isomorphism of $q$-inverting
pairs from
 $K,K^*$ to
$K', K^{*\prime}$.
\end{definition}

\noindent
The main results of this paper 
are contained
in the following two theorems and subsequent
remark.

\begin{theorem}
\label{thm:easy}
Let $V$ denote a finite-dimensional irreducible
$\boxtimes_q$-module of type 1. Then the generators
$x_{02}, x_{13}$ act on $V$ as a $q$-inverting pair.
\end{theorem}

\begin{theorem}
\label{thm:hard}
Let $V$ denote a vector space over $\K$ with
finite positive dimension and let $K,K^*$ denote
a $q$-inverting pair on $V$.
Then there exists a unique $\boxtimes_q$-module
structure on $V$ such that $x_{02}, x_{13}$ act
on $V$ as $K,K^*$ respectively.
This module structure is irreducible and type 1.
\end{theorem}

\begin{remark} 
\rm 
Combining Theorem \ref{thm:easy} and
Theorem
\ref{thm:hard} we get a bijection between the
following two sets:
\begin{enumerate}
\item[{\rm (i)}]  
the isomorphism classes of finite-dimensional irreducible
$\boxtimes_q$-modules of type 1;
\item[{\rm (ii)}]  
the isomorphism classes of $q$-inverting pairs.
\end{enumerate}
\end{remark}


\noindent The proof of 
Theorem
\ref{thm:easy} and 
 Theorem
\ref{thm:hard} will take up Sections 5--8.

\section{The $\Z_4$ action}

\noindent
In this section we display an
action of the group $\Z_4$ on the set of $q$-inverting pairs.

\medskip
\noindent Referring to the $q$-inverting pair
$K,K^*$ on $V$ from
 Definition
\ref{def:inv}, if we replace
\begin{eqnarray*}
K; K^*; \lbrace V_i\rbrace_{i=0}^d; \lbrace V^*_i\rbrace_{i=0}^{\delta}
\end{eqnarray*}
by 
\begin{eqnarray*}
K^*; K^{-1}; 
\lbrace V^*_i\rbrace_{i=0}^{\delta}; 
\lbrace V_{d-i}\rbrace_{i=0}^d
\end{eqnarray*}
then the axioms in Definition \ref{def:inv}(i)--(iii) 
still hold; therefore  
the pair $K^*, K^{-1}$ is a $q$-inverting pair on $V$.
Consider the map $\varrho$ which takes each
$q$-inverting pair $K,K^*$ to the $q$-inverting pair
 $K^*, K^{-1}$.
The  map $\varrho$ is a permutation on the set of $q$-inverting pairs,
and $\varrho^4=1$. Therefore $\varrho$ induces an action of $\Z_4$
on the set of $q$-inverting pairs.
We record a result for later use.

\begin{corollary}
\label{cor:four}
Let $K,K^*$ denote a $q$-inverting pair.
Then each of the following is a $q$-inverting pair:
\begin{eqnarray}
\label{eq:four}
K,K^*;
\qquad \qquad 
K^*,K^{-1};
\qquad \qquad 
K^{-1},K^{*-1};
\qquad \qquad 
K^{*-1},K.
\end{eqnarray}
\end{corollary}
\noindent {\it Proof:} Repeatedly apply $\varrho$
to the $q$-inverting pair $K,K^*$.
\hfill $\Box $ \\ 

\begin{remark}
\rm
The 
 $q$-inverting pairs 
(\ref{eq:four}) 
might not be mutually nonisomorphic.
\end{remark}

\section{Some linear algebra}

\noindent In this section we obtain
some linear algebraic results that we
will need to prove Theorem \ref{thm:easy}
and Theorem \ref{thm:hard}.
We will use the following concepts.
Let $V$ denote a
vector space over $\K$ with finite positive dimension 
and let $A:V\to V$ denote a linear transformation.
For $\theta \in \K$ we define
\begin{eqnarray*}
V_A(\theta) = \lbrace v \in V \,|\,Av = \theta v\rbrace.
\end{eqnarray*}
Observe that $\theta$ is an eigenvalue of $A$ if and only
if $V_A(\theta)\not=0$, and in this case
$V_A(\theta)$ is the corresponding eigenspace.
The sum
$\sum_{\theta \in \K} V_A(\theta)$ is direct.
Moreover 
this sum is equal to $V$ 
if and only if
$A$ is semisimple.

\begin{lemma}
\label{lem:key}
Let $V$ denote a 
vector space over $\K$ with finite positive dimension.
Let $A:V\to V$ and $B:V\to V$ denote linear transformations.
Then for all nonzero $\theta \in \K$ the following 
are equivalent:
\begin{enumerate}
\item
The expression
$
A^3B-\lbrack 3\rbrack_q A^2BA
+\lbrack 3\rbrack_q ABA^2
-BA^3
$
vanishes on $V_A(\theta)$.
\item
$B V_A(\theta) \subseteq 
 V_A(q^2\theta)
+
 V_A(\theta)
+
V_A(q^{-2}\theta)
$.
\end{enumerate}
\end{lemma}
\noindent {\it Proof:} 
For $v \in V_A(\theta)$ we
have 
\begin{eqnarray*}
&&(A^3B-\lbrack 3\rbrack_q A^2BA
+\lbrack 3\rbrack_q ABA^2
-BA^3)v
\\
&& \qquad = 
(A^3-\theta \lbrack 3\rbrack_q A^2
+\theta^2\lbrack 3\rbrack_q A
-\theta^3I)Bv
\qquad \qquad \mbox{since}\;Av=\theta v
\\
&&\qquad =
(A-q^2 \theta I)
(A- \theta I)
(A-q^{-2} \theta I)Bv,
\end{eqnarray*}
where $I:V\to V$ is the identity map.
The scalars $q^{2}\theta, \theta, q^{-2}\theta$
are mutually distinct since $\theta\not=0$ and
since $q$ is not a root of $1$.
The result follows.
\hfill $\Box $

\begin{lemma}
\label{lem:qweyl}
Let $V$ denote a vector space over $\K$ with finite positive
dimension. Let $A:V\to V$ and $B:V\to V$ denote linear transformations.
Then for all nonzero $\theta \in \K$ the following are equivalent:
\begin{enumerate}
\item
The expression
$
qAB-q^{-1}BA-(q-q^{-1})I
$
vanishes on $V_A(\theta)$.
\item
$(B -\theta^{-1}I)V_A(\theta) \subseteq 
 V_A(q^{-2}\theta)
$.
\end{enumerate}
\end{lemma}
\noindent {\it Proof:} 
For $v \in V_A(\theta)$ we
have 
\begin{eqnarray*}
(qAB-q^{-1}BA-(q-q^{-1})I)v=
q(A-q^{-2}\theta I)(B-\theta^{-1}I)v
\end{eqnarray*}
and the result follows.
\hfill $\Box $ \\

\begin{lemma}
\label{lem:qweyl2}
Let $V$ denote a vector space over $\K$ with finite positive
dimension. Let $A:V\to V$ and $B:V\to V$ denote linear transformations.
Then for all nonzero $\theta \in \K$ the following are equivalent:
\begin{enumerate}
\item
The expression
$
qAB-q^{-1}BA-(q-q^{-1})I
$
vanishes on $V_B(\theta)$.
\item
$(A -\theta^{-1}I)V_B(\theta) \subseteq 
 V_B(q^{2}\theta)
$.
\end{enumerate}
\end{lemma}
\noindent {\it Proof:} 
In Lemma
\ref{lem:qweyl} replace
$(A,B,q)$ by $(B,A,q^{-1})$.
\hfill $\Box $  \\

\section{From $\boxtimes_q$-modules to $q$-inverting pairs}

Our goal in this section is to prove
Theorem
\ref{thm:easy}. We start with some comments on
$\boxtimes_q$-modules.

\begin{definition}
\label{def:decij}
\rm
Let $V$ denote a
finite-dimensional irreducible $\boxtimes_q$-module of type
$1$ and diameter $d$. For each generator $x_{rs}$ of $\boxtimes_q$
we define a decomposition of $V$ which we call $\lbrack r,s\rbrack$.
The decomposition $\lbrack r,s\rbrack $ has diameter $d$.
For $0 \leq i\leq d$ the $i^{th}$ component of
$\lbrack r,s\rbrack$ is the eigenspace of $x_{rs}$ on
$V$ associated with the eigenvalue $q^{d-2i}$. 
\end{definition}

\begin{note}\rm
With reference to Definition
\ref{def:decij}, for $r \in \Z_4$ the decomposition
$\lbrack r,r+2\rbrack $ is the inversion of
the decomposition $\lbrack r+2,r\rbrack$.
\end{note}

\begin{proposition}
\cite[Proposition 13.3]{qtet}
\label{lem:shape}
Let $V$ denote a finite-dimensional irreducible $\boxtimes_q$-module
of type 1 and diameter $d$. Choose a generator $x_{rs}$ of $\boxtimes_q$
and consider the corresponding decomposition $\lbrack r,s\rbrack$ of $V$  
from Definition
\ref{def:decij}.
Then the shape of this decomposition is independent
of the choice of generator. Denoting the shape by
$\lbrace \rho_i\rbrace_{i=0}^d$ we
have $\rho_i=\rho_{d-i}$ for $0 \leq i\leq d$.
\end{proposition}

\begin{definition} \rm
\label{def:shapev}
Let $V$ denote a finite-dimensional irreducible $\boxtimes_q$-module
of type 1 and diameter $d$. By the {\it shape} of $V$
we mean the sequence 
$\lbrace \rho_i\rbrace_{i=0}^d$
from
Proposition
\ref{lem:shape}.
\end{definition}

\begin{theorem}
\cite[Theorem 14.1]{qtet}
\label{thm:sixdecp}
Let $V$ denote a finite-dimensional irreducible
$\boxtimes_q$-module of type 1 and diameter $d$.
Let 
$\lbrace U_i\rbrace_{i=0}^d$ denote a decomposition of
$V$ from Definition
\ref{def:decij}. Then for $r \in \Z_4$ and
for $0 \leq i\leq d$ the action of $x_{r,r+1}$ on
$U_i$ is given as follows.
\medskip

\centerline{
\begin{tabular}[t]{c|c}
       {\rm decomposition} & {\rm action of $x_{r,r+1}$ on $U_i$}
 \\ \hline  \hline
	$\lbrack r,r+1\rbrack$ & $(x_{r,r+1}-q^{d-2i}I)U_i=0$    
	\\
	$\lbrack r+1,r+2\rbrack$ & 
              $(x_{r,r+1}-q^{2i-d}I)U_i \subseteq U_{i-1}$   \\
	$\lbrack r+2,r+3\rbrack$ &
	$x_{r,r+1}U_i\subseteq U_{i-1}+U_i+U_{i+1}$  \\ 
	$\lbrack r+3,r\rbrack$ & 
	$(x_{r,r+1}-q^{2i-d}I)U_i\subseteq U_{i+1}$
	\\
        $\lbrack r,r+2\rbrack$ &
	$
	(x_{r,r+1}-q^{d-2i}I)U_i\subseteq U_{i-1}$
	\\ 
	$ \lbrack r+1,r+3\rbrack $ & 
	$
	(x_{r,r+1}-q^{2i-d}I)U_i\subseteq U_{i-1}$
	\end{tabular}}
\medskip
\noindent 
\end{theorem}

\begin{theorem}
\cite[Theorem 14.2]{qtet}
\label{thm:sixdecp2}
Let $V$ denote a finite-dimensional irreducible
$\boxtimes_q$-module of type 1 and diameter $d$.
Let 
$\lbrace U_i\rbrace_{i=0}^d$
denote a decomposition of
$V$ from Definition
\ref{def:decij}. Then for $r \in \Z_4$ and
for $0 \leq i\leq d$ the action of $x_{r,r+2}$ on
$U_i$ is given as follows.
\medskip

\centerline{
\begin{tabular}[t]{c|c}
       {\rm decomposition} & {\rm action of $x_{r,r+2}$ on $U_i$}
 \\ \hline  \hline
	$\lbrack r,r+1\rbrack$ &
$(x_{r,r+2}-q^{d-2i}I)U_i \subseteq U_0+\cdots+U_{i-1}$ 
	\\
	$\lbrack r+1,r+2\rbrack$ & 
       $(x_{r,r+2}-q^{d-2i}I)U_i \subseteq U_{i+1}+\cdots+U_{d}$     \\
	$\lbrack r+2,r+3\rbrack$ &
	$(x_{r,r+2}-q^{2i-d}I)U_i\subseteq U_{i-1}$
         \\ 
	$\lbrack r+3,r\rbrack$ & 
	$(x_{r,r+2}-q^{2i-d}I)U_i\subseteq U_{i+1}$
	\\ 
        $\lbrack r,r+2\rbrack$ &
	$
	(x_{r,r+2}-q^{d-2i}I)U_i=0$
	\\ 
	$ \lbrack r+1,r+3\rbrack $ & 
	$
	x_{r,r+2}U_i\subseteq U_{i-1}+\cdots+U_d$
	\end{tabular}}
\medskip
\noindent 
\end{theorem}

\noindent We recall the notion of a {\it flag}.
Let $V$ denote a vector space over $\K$ with finite
positive dimension. Let $\lbrace s_i\rbrace_{i=0}^d$ denote a sequence
of positive integers whose sum is the dimension of $V$.
By a {\it flag on $V$ of shape $\lbrace s_i\rbrace_{i=0}^d$} we mean
a nested sequence $F_0\subseteq F_1\subseteq \cdots \subseteq
F_d$ of subspaces of $V$ such that the dimension of $F_i$
is $s_0+\cdots+s_i$ for $0 \leq i \leq d$. We call
$F_i$ the {\it $i^{th}$ component} of the flag. We call
$d$ the {\it diameter} of the flag. We observe 
$F_d=V$.

\medskip
\noindent The following construction yields a flag on $V$.
Let $\lbrace U_i\rbrace_{i=0}^d$ denote a decomposition of $V$
of shape $\lbrace s_i\rbrace_{i=0}^d$. Define
\begin{eqnarray*}
F_i=U_0+U_1+\cdots+U_i \qquad \qquad (0 \leq i\leq d).
\end{eqnarray*}
Then the sequence $F_0 \subseteq F_1 \subseteq \cdots \subseteq F_d$
is a flag on $V$ of shape $\lbrace s_i\rbrace_{i=0}^d$.
We say this flag is {\it induced} by the decomposition
 $\lbrace U_i\rbrace_{i=0}^d$.

\medskip
\noindent We now recall what it means for two flags to be {\it opposite}.
Suppose we are given two flags on $V$ with the same diameter: 
 $F_0 \subseteq F_1 \subseteq \cdots \subseteq F_d$
and 
 $F'_0 \subseteq F'_1 \subseteq \cdots \subseteq F'_d$.
We say that these flags are {\it opposite} whenever there exists
a decomposition $\lbrace U_i\rbrace_{i=0}^d$ of $V$ such that
\begin{eqnarray*}
F_i=U_0+U_1+\cdots+U_i, \qquad \qquad 
F'_i=U_d+U_{d-1}+\cdots+U_{d-i}
\end{eqnarray*}
for $0 \leq i\leq d$. In this case
\begin{eqnarray}
\label{eq:zero}
F_i\cap F'_j = 0 \qquad \mbox{if}\quad i+j<d \qquad \qquad (0 \leq i,j\leq d)
\end{eqnarray}
and
\begin{eqnarray}
\label{eq:recover}
U_i=F_i\cap F'_{d-i} \qquad \qquad (0 \leq i\leq d).
\end{eqnarray}
In particular the decomposition $\lbrace U_i\rbrace_{i=0}^d$ is uniquely
determined by the given flags.

\medskip
\noindent We  now return our attention to
$\boxtimes_q$.

\begin{theorem}
\cite[Theorem 16.1]{qtet}
\label{thm:fourflags}
Let $V$ denote a finite-dimensional irreducible $\boxtimes_q$-module
of type 1 and diameter $d$. Then there exists a collection of
flags on $V$, denoted $\lbrack n \rbrack, n \in \Z_4$, that have 
the following property:
for each generator $x_{rs}$ of $\boxtimes_q$ the decomposition
$\lbrack r,s\rbrack$ of $V$ induces $\lbrack r\rbrack$ and
the inversion of $\lbrack r,s\rbrack$ induces $\lbrack s\rbrack$.
\end{theorem}

\begin{lemma}
\cite[Lemma 16.2]{qtet}
Let $V$ denote a finite-dimensional irreducible
$\boxtimes_q$-module of type 1. Then for $n \in \Z_4$ the shape
of the flag $\lbrack n \rbrack $ coincides with the shape of $V$.
\end{lemma}

\begin{theorem} 
\cite[Theorem 16.3]{qtet}
Let $V$ denote a finite-dimensional 
irreducible $\boxtimes_q$-module of
type 1. Then the flags $\lbrack n \rbrack$, $n \in \Z_4$ on $V$  
from Theorem \ref{thm:fourflags}
are mutually opposite.
\end{theorem}

\begin{theorem}
\cite[Theorem 16.4]{qtet}
\label{thm:flagdec}
Let $V$ denote a finite-dimensional irreducible $\boxtimes_q$-module
of type 1 and diameter $d$. Pick a generator $x_{rs}$ of $\boxtimes_q$
and consider the corresponding decomposition $\lbrack r,s\rbrack $
of $V$ from Definition \ref{def:decij}.
For $0 \leq i\leq d$ the $i^{th}$ component of $\lbrack r,s\rbrack$
is the intersection of the following two sets:
\begin{enumerate}
\item[{\rm (i)}]  
component $i$ of the flag $\lbrack r\rbrack$;
\item[{\rm (ii)}] 
component $d-i$ of the flag $\lbrack s\rbrack$.
\end{enumerate}
\end{theorem}

\begin{proposition}
\label{prop:w}
Let $V$ denote a finite-dimensional irreducible $\boxtimes_q$-module
of type 1. Let $W$ denote a nonzero subspace of $V$ such
that $x_{02}W\subseteq W$ and
$x_{13}W\subseteq W$. Then $W=V$.
\end{proposition}
\noindent {\it Proof:}
Without loss we may assume that $W$ is irreducible as a module
for $x_{02},x_{13}$.
Since $x_{02}x_{20}=1$ and
$x_{13}x_{31}=1$ 
we find
$W$ is invariant under each of $x_{20},x_{31}$.
Therefore 
$W$ is invariant under $x_{r,r+2}$ for $r \in \Z_4$.
For the moment fix $r \in \Z_4$
and let 
 $\lbrace U_i\rbrace_{i=0}^d$ denote the decomposition
$\lbrack r,r+2\rbrack$.
Recall that $x_{r,r+2}$
is semisimple on $V$ with eigenspaces
$U_0,U_1,\ldots,U_d$. 
By this and since $W$ is invariant under
 $x_{r,r+2}$ we find
\begin{eqnarray*}
\label{eq:wsum}
W=\sum_{i=0}^d W\cap U_i.
\end{eqnarray*}
Since $W\not=0$  there exists an integer $i$ $(0 \leq i \leq d)$
such that $W\cap U_i\not=0$.
Define
\begin{eqnarray*}
\label{eq:m}
m_r= \mbox{\rm min}\lbrace i \;|\;0 \leq i\leq d, 
\quad W\cap U_i\not=0\rbrace.
\end{eqnarray*}
We claim that $m_r$ is independent of $r$ for $r \in \Z_4$.
Suppose the claim is false.
Then there exists $r \in \Z_4$ such that
 $m_r>m_{r+1}$. By  construction
the space $W$ is contained
 in component $d-m_r$ of the flag $\lbrack r+2\rbrack$.
By  construction $W$ has nonzero intersection
with component $m_{r+1}$
of the flag $\lbrack r+1\rbrack$. Since $m_r>m_{r+1}$
the component $d-m_r$ of
$\lbrack r+2\rbrack$ has zero intersection with  component
$m_{r+1}$ of $\lbrack r+1\rbrack$, for a contradiction.
We have proved the claim.
For the rest of this proof let $m$ denote the common value
of $m_r$ for  $r\in \Z_4$.
The claim implies that for all $r \in \Z_4$ the component
$d-m$ of the flag $\lbrack r \rbrack $ contains $W$,
and component $m$ of $\lbrack r \rbrack$ has nonzero intersection
with $W$.
We can now easily show  $W=V$.
Since the $\boxtimes_q$-module 
$V$ is irreducible and $W\not=0$ it suffices
to show that $W$ is invariant under $\boxtimes_q$. 
We mentioned earlier that 
$W$ is invariant under $x_{r,r+2}$ for $r \in \Z_4$.
We now show that $W$ is invariant under
$x_{r,r+1}$ for $r \in \Z_4$.
Let $r$ be given and 
let $W'$ denote the span of the set of vectors in 
$W$ that are eigenvectors for $x_{r,r+1}$. By the construction
$W'\subseteq W$ and $x_{r,r+1}W'\subseteq W'$. We show
$W'=W$. To this end we show that
$W'$ is nonzero and invariant under each of $x_{02}, x_{13}$.
We now show  $W'\not=0$. By the comment after the preliminary 
claim, $W$ has nonzero intersection with component $m$
of the flag $\lbrack r \rbrack $ and $W$ is contained in
component $d-m$ of the flag $\lbrack r+1 \rbrack$. By 
Theorem
\ref{thm:flagdec} the intersection of component
$m$ of $\lbrack r\rbrack $ and component $d-m$ of
$\lbrack r+1\rbrack$ is equal to component $m$ of the
decomposition $\lbrack r,r+1\rbrack$, which is an eigenspace
for $x_{r,r+1}$. The intersection of $W$ with this eigenspace
is nonzero and contained in $W'$, so $W'\not=0$.
We now show that 
$W'$ is invariant under each of $x_{02}, x_{13}$.
Since $x_{02}x_{20}=1$ and $x_{13}x_{31}=1$
it suffices to show that $W'$ is invariant under
 $x_{r+1,r+3}$ and $x_{r+2,r}$.
We now show that 
$W'$ is invariant under $x_{r+1,r+3}$.
To this end we pick $v \in W'$ and show $x_{r+1,r+3}v \in W'$.
Without loss we may assume that $v$ is an eigenvector for $x_{r,r+1}$;
let $\theta$ denote the corresponding eigenvalue. Then
$\theta\not=0$ by the comment
at the end of Section 2.
Since $v\in W'$ and
 $W'\subseteq W$
we have  $v \in W$. The space $W$ is invariant under
$x_{r+1,r+3}$ so
$x_{r+1,r+3}v \in W$.
By these comments $(x_{r+1,r+3}-\theta^{-1}I)v \in W$. By Lemma
\ref{lem:qweyl} (with $A=x_{r,r+1}$ and $B=x_{r+1,r+3}$) the vector
$(x_{r+1,r+3}-\theta^{-1}I)v$ is contained in an eigenspace
of $x_{r,r+1}$ 
so $(x_{r+1,r+3}-\theta^{-1}I)v\in W'$. By this and since
$v \in W'$ we have $x_{r+1,r+3}v \in W'$.
We have now shown that $W'$ is invariant under $x_{r+1,r+3}$ as desired.
Next we show that $W'$ is invariant under $x_{r+2,r}$.
To this end we pick $u \in W'$ and show $x_{r+2,r}u \in W'$.
Without loss we may assume that $u$ is an eigenvector for $x_{r,r+1}$;
let $\eta$ denote the corresponding eigenvalue. Then
$\eta\not=0$ by
the comment at the end of Section 2.
Recall 
$u\in W'$
and $W'\subseteq W$ 
so $u \in W$. The space $W$ is invariant under
$x_{r+2,r}$ so $x_{r+2,r}u \in W$.
By these comments $(x_{r+2,r}-\eta^{-1}I)u \in W$. By
 Lemma \ref{lem:qweyl2} (with $A=x_{r+2,r}$, $B=x_{r,r+1}$, $\theta=\eta$) 
 the vector
$(x_{r+2,r}-\eta^{-1}I)u$ is contained in an eigenspace of $x_{r,r+1}$
so $(x_{r+2,r}-\eta^{-1}I)u\in W'$. By this and since
$u \in W'$ we have $x_{r+2,r}u \in W'$.
We have now shown that $W'$ is invariant under $x_{r+2,r}$ as desired.
From our above comments
$W'$ is nonzero and invariant under each of
$x_{02},x_{13}$. Now $W'=W$ by the irreducibility of $W$,
so $x_{r,r+1}W\subseteq W$. 
We have now shown that $W$ is invariant under $x_{r,r+1}$ and
$x_{r,r+2}$ for $r \in \Z_4$.
Therefore 
$W$ is $\boxtimes_q$-invariant. The $\boxtimes_q$-module
$V$ is
irreducible and $W\not=0$ so  $W=V$.
\hfill $\Box $ \\

\noindent It is now a simple matter to prove Theorem
\ref{thm:easy}. 

\medskip
\noindent {\it Proof of Theorem 
\ref{thm:easy}}:
Define the linear transformation  $K:V\to V$
(resp. $K^*:V\to V$)
to be the action of $x_{02}$ (resp. $x_{13}$) on $V$.
We show that $K,K^*$ is a $q$-inverting pair on $V$.
To do this we show that
$K,K^*$ satisfy the conditions (i)--(iii) of
Definition \ref{def:inv}.
Concerning Definition \ref{def:inv}(i),
we denote the decomposition $\lbrack 0,2\rbrack$ 
by $\lbrace V_i \rbrace_{i=0}^d$
 and show that this decomposition satisfies
(\ref{eq:k1})--(\ref{eq:k3}).
Line (\ref{eq:k1}) 
is satisfied by the construction.
To get 
 (\ref{eq:k2}),  
 (\ref{eq:k3})  
we refer to the last row 
in the table of Theorem
\ref{thm:sixdecp2}.
Line (\ref{eq:k2}) holds
 by
that row
 (with $r=1$) and since the decomposition
$\lbrack 2,0 \rbrack$ is the inversion of 
$\lbrack 0,2 \rbrack$.
Line 
(\ref{eq:k3}) holds by
that row
 (with $r=3$) and since $x_{13}x_{31}=1$.
We have now shown that
$K,K^*$ satisfy 
Definition \ref{def:inv}(i).
Concerning Definition \ref{def:inv}(ii),
we denote the decomposition $\lbrack 1,3\rbrack$ 
by $\lbrace V^*_i \rbrace_{i=0}^d$
 and show that this decomposition satisfies
(\ref{eq:ks1})--(\ref{eq:ks3}).
Line (\ref{eq:ks1}) holds
by the construction.
To get 
 (\ref{eq:ks2}),  
 (\ref{eq:ks3})  
 we refer  to the last row 
in the table of Theorem
\ref{thm:sixdecp2}.
Line (\ref{eq:ks2}) holds by
that row
 (with $r=0$).
Line 
(\ref{eq:ks3}) holds
by
that row
 (with $r=2$), since 
$x_{02}x_{20}=1$
 and since 
 the decomposition
$\lbrack 3,1 \rbrack$ is the inversion of 
$\lbrack 1,3 \rbrack$.
We have now shown that
$K,K^*$ satisfy 
Definition \ref{def:inv}(ii).
The maps $K,K^*$ satisify Definition \ref{def:inv}(iii)
by Proposition
\ref{prop:w}.
We have now verified  that
$K,K^*$ satisfy 
Definition \ref{def:inv}(i)--(iii) so $K,K^*$ is a $q$-inverting pair on $V$.
The result follows.
\hfill $\Box $ \\

\section{From $q$-inverting pairs to $\boxtimes_q$-modules}

Our goal in this section is to  prove Theorem
\ref{thm:hard}. On our way to this goal
we will show that the integers $d$ and $\delta$ from
Definition \ref{def:inv} are equal.

\begin{definition}
\label{def:vij}
\rm
With reference to 
Definition
\ref{def:inv} 
we set 
\begin{eqnarray}
V_{ij} =  
(V_0+\cdots + V_i)\cap (V^*_0+\cdots + V^*_j)
\label{eq:vij}
\end{eqnarray}
for all integers $i,j$. We interpret the sum on the left in
(\ref{eq:vij}) to be 0 (resp. $V$) if $i<0$  (resp. $i>d$).
We interpret the sum on the right in
(\ref{eq:vij}) to be 0 (resp. $V$) if $j<0$  (resp. $j>\delta$).
\end{definition}

\begin{lemma}
\label{lem:basic}
With reference to 
Definition
\ref{def:inv} and 
Definition \ref{def:vij}, the following (i), (ii) hold.
\begin{enumerate}
\item
$V_{i\delta} = 
 V_0+\cdots+V_i \qquad  (0 \leq i \leq d)$.
\item
$V_{dj} = 
 V^*_0+\cdots+V^*_j\qquad (0 \leq j\leq \delta)$.
\end{enumerate}
\end{lemma}
\noindent {\it Proof:} 
(i) Set $j=\delta$ in 
(\ref{eq:vij}) and recall
$V=V^*_0+\cdots + V^*_\delta$.
\\
(ii) Set $i=d$ in 
(\ref{eq:vij}) and use
$V=V_0+\cdots +V_d$.
\hfill $\Box $

\begin{lemma}
\label{lem:kact}
With reference to 
Definition
\ref{def:inv} 
and Definition \ref{def:vij}, the
following (i), (ii) hold
for $0 \leq i\leq d$ and $0 \leq j\leq \delta$.
\begin{enumerate}
\item
$(K^{-1}-q^{2i-d}I)V_{ij} \subseteq V_{i-1,j+1}$.
\item
$(K^*-q^{\delta-2j}I)V_{ij} \subseteq V_{i+1,j-1}$.
\end{enumerate}
\end{lemma}
\noindent {\it Proof:} 
(i) Using 
(\ref{eq:k1})
we find
\begin{eqnarray}
(K^{-1}-q^{2i-d}I)\sum_{h=0}^i V_h =
\sum_{h=0}^{i-1} V_h.
\label{eq:kinv1}
\end{eqnarray}
Using
(\ref{eq:ks3})  we find
\begin{eqnarray}
(K^{-1}-q^{2i-d}I)\sum_{h=0}^j V^*_h \subseteq 
\sum_{h=0}^{j+1} V^*_h.
\label{eq:kinv2}
\end{eqnarray}
Evaluating 
$(K^{-1}-q^{2i-d}I)V_{ij}$
using 
(\ref{eq:vij})--(\ref{eq:kinv2})
we find it is contained in 
$V_{i-1,j+1}$.
\\
\noindent (ii)
 Using 
(\ref{eq:k2})
we find
\begin{eqnarray}
(K^*-q^{\delta-2j}I)\sum_{h=0}^i V_h \subseteq 
\sum_{h=0}^{i+1} V_h.
\label{eq:ksinv1}
\end{eqnarray}
Using
(\ref{eq:ks1})  we find
\begin{eqnarray}
(K^*-q^{\delta -2j}I)\sum_{h=0}^j V^*_h 
=\sum_{h=0}^{j-1} V^*_h.
\label{eq:ksinv2}
\end{eqnarray}
Evaluating 
$(K^*-q^{\delta-2j}I)V_{ij}$
using 
(\ref{eq:vij}) and (\ref{eq:ksinv1}),
(\ref{eq:ksinv2}) 
we find it is contained in 
$V_{i+1,j-1}$.
\hfill $\Box $  

\begin{lemma}
\label{lem:dvsdel}
The scalars 
 $d$ and $\delta $ 
from Definition
\ref{def:inv}
 are equal. Moreover,
with reference to 
Definition \ref{def:vij}, 
\begin{eqnarray}
V_{ij}= 0 \quad \mbox{if}\quad i+j<d\qquad  \qquad (0 \leq i,j\leq d).
\label{eq:ijd}
\end{eqnarray}
\end{lemma}
\noindent {\it Proof:} 
For all nonnegative integers $r$ such that
$r\leq d$ and $r \leq \delta$ we define
\begin{eqnarray}
W_r=V_{0r}+V_{1,r-1}+\cdots +V_{r0}.
\label{eq:wr}
\end{eqnarray}
We have $K^{-1}W_r \subseteq W_r$
by Lemma
\ref{lem:kact}(i) so
$KW_r \subseteq W_r$. We have
$K^*W_r \subseteq W_r$ by Lemma
\ref{lem:kact}(ii).
Now $W_r=0$ or $W_r=V$ in view of Definition
\ref{def:inv}(iii).
Suppose for the moment that $r\leq d-1$.
Each term on the right in 
(\ref{eq:wr}) is contained in
$V_0+\cdots + V_r$ so
$W_r \subseteq 
V_0+\cdots + V_r$.
Therefore
$W_r\not=V$ so $W_r=0$.
Next suppose $r=d$.
Then 
$V_{d0}\subseteq W_r$.
 Recall $V_{d0}=V^*_0$ by Lemma
\ref{lem:basic}(ii) and $V^*_0\not=0$ so
 $V_{d0}\not=0$.
Now $W_r\not=0$ so $W_r=V$.
We have now shown that $W_r=0$ if $r\leq d-1$ and
$W_r=V$ if $r=d$.
Similarly
$W_r=0$ if $r\leq \delta-1$ and
$W_r=V$ if $r=\delta$.
Now $d=\delta$; otherwise 
we take $r=\mbox{min}(d,\delta)$ in our above comments
and find $W_r$ is both $0$ and $V$, for a contradiction.
The result follows.
\hfill $\Box $

\begin{definition}
\label{def:ui}
With reference to 
Definition
\ref{def:inv}, for $0 \leq i \leq d$ we define
\begin{eqnarray*}
U_i = (V_0+\cdots + V_i)\cap (V^*_0+\cdots + V^*_{d-i}).
\end{eqnarray*}
We observe $U_i$ is equal to the space 
$V_{i,d-i}$ from Definition \ref{def:vij}.
\end{definition}

\begin{lemma}
\label{lem:udec}
With reference to
Definition \ref{def:inv} and
Definition \ref{def:ui}, the sequence
$\lbrace U_i\rbrace_{i=0}^d$ is a decomposition of $V$.
\end{lemma}
\noindent {\it Proof:} 
We first show 
\begin{eqnarray}
V= \sum_{i=0}^d U_i.
\label{eq:usumv}
\end{eqnarray}
Let $W$ denote the sum on the right
in 
(\ref{eq:usumv}).
We have $K^{-1}W\subseteq W$ by 
Lemma \ref{lem:kact}(i) so $KW\subseteq W$.
We have $K^*W\subseteq W$ by 
Lemma \ref{lem:kact}(ii).
Now 
$W=0$ or $W=V$ in view of
Definition
\ref{def:inv}(iii).
The space
$W$ contains 
 $U_{0}$ and $U_{0}=V_0$ is nonzero
 so
$W\not=0$.
Therefore $W=V$ and
(\ref{eq:usumv}) follows.
Next we show that the sum 
(\ref{eq:usumv}) is direct.  To do this,
we show
\begin{eqnarray*}
 (U_{0}+U_{1}+\cdots +U_{i-1})\cap U_{i} =0
\end{eqnarray*}
for $1 \leq i \leq d$.
Let  $i$ be given. From the construction 
\begin{eqnarray*}
U_{j} \subseteq V_0+V_1+\cdots +V_{i-1}
\end{eqnarray*}
for $0 \leq j \leq i-1$, and
\begin{eqnarray*}
U_{i}\subseteq V^*_0+V^*_1+\cdots +V^*_{d-i}.
\end{eqnarray*}
Therefore 
\begin{eqnarray*}
 &&(U_{0}+U_{1}+\cdots +U_{i-1})\cap U_{i}
\\
 && \qquad  \subseteq
 (V_0+V_1+\cdots +V_{i-1})\cap
(V^*_0+V^*_1+\cdots +V^*_{d-i}) \qquad \qquad 
\\
&& \qquad  = V_{i-1,d-i}
\\
&& \qquad  = 0
\end{eqnarray*}
in view of
(\ref{eq:ijd}).
We have now shown
that the sum (\ref{eq:usumv}) is direct.
Next we show that each of
$U_0,\ldots, U_d$ is nonzero.
Suppose there exists an integer $i$ $(0 \leq i \leq d)$
such that $U_i=0$. Observe that $i\not=0$ since
$U_0=V_0$ is nonzero, and $i\not=d$ since $U_d=V^*_0$
is nonzero. Set
\begin{eqnarray*}
U= U_0+U_1 + \cdots+U_{i-1}
\end{eqnarray*}
and observe $U\not=0$ and $U\not=V$ by our above
remarks. By Lemma
\ref{lem:kact}(i)
we find $K^{-1}U\subseteq U$ 
so $KU \subseteq U$.
By 
 Lemma
\ref{lem:kact}(ii)
and since $U_i=0$ 
we find $K^*U\subseteq U$.
Now $U=0$ or $U=V$ by Definition \ref{def:inv}(iii),
for a contradiction.
We conclude that 
each of $U_0,\ldots, U_d$ is nonzero. 
We have now shown that the sequence $\lbrace U_i\rbrace_{i=0}^d$
is
a decomposition of $V$.
\hfill $\Box $  

\begin{definition}
\rm
We call the decomposition 
$\lbrace U_i\rbrace_{i=0}^d$ 
from Lemma
\ref{lem:udec}
the {\it split decomposition} of $V$
associated with $K,K^*$.
\end{definition}

\begin{lemma}
\label{lem:kks}
With reference to Definition
\ref{def:inv} and Definition
\ref{def:ui}, the following (i), (ii) hold
for $0 \leq i \leq d$.
\begin{enumerate}
\item
$(K^{-1}-q^{2i-d}I)U_i\subseteq U_{i-1}$.
\item
$(K^*-q^{2i-d}I)U_i\subseteq U_{i+1}$.
\end{enumerate}
\end{lemma}
\noindent {\it Proof:} 
Immediate from Lemma
\ref{lem:kact}
 and Definition
\ref{def:ui}.
\hfill $\Box $  

\begin{lemma}
\label{lem:uv}
With reference to Definition
\ref{def:inv} and Definition
\ref{def:ui}, the following (i), (ii) hold
for $0 \leq i \leq d$.
\begin{enumerate}
\item
$U_0+\cdots+U_i=V_0+\cdots+V_i$.
\item
$U_{d-i}+\cdots+U_d=V^*_0+\cdots+V^*_i$.
\end{enumerate}
\end{lemma}
\noindent {\it Proof:}
\noindent (i)
Let $X_i=\sum_{j=0}^i U_j$
and
$X'_i=\sum_{j=0}^i V_j$.
We show $X_i=X'_i$.
Define $R_i = \prod_{j=0}^i (K^{-1}-q^{2j-d}I)$.
Then $X'_i = \lbrace v \in V\,|\,R_iv=0\rbrace$,
and $R_iX_i=0$ by
Lemma \ref{lem:kks}(i),
so $X_i \subseteq X'_i$.
Now define $T_i= \prod_{j=i+1}^{d}(K^{-1}-q^{2j-d}I)$.
Observe that $T_iV=X'_i$, and
 $T_iV\subseteq X_i$ by
 Lemma \ref{lem:kks}(i),
 so $X'_i \subseteq X_i$.
 By these comments $X_i=X'_i$.
 \\
 \noindent (ii) Similar to the proof of (i) above.
 \hfill $\Box $ \\

\begin{definition}
\label{def:spltop}
\rm
With reference to Definition
\ref{def:inv}, by the {\it split operator}
for $K,K^*$ we mean the linear transformation
$S: V\to V$
such that for $0 \leq i \leq d$ the space
$U_i$ from Definition
\ref{def:ui}
is the eigenspace for $S$ with eigenvalue
$q^{d-2i}$.
\end{definition}

\noindent We are now ready to prove
Theorem \ref{thm:hard}.

\medskip
\noindent {\it Proof of Theorem 
 \ref{thm:hard}:}
We let the generators $x_{rs}$ of 
$\boxtimes_q$
act on 
$V$ as follows.
We let $x_{02}$ (resp. $x_{13}$)
(resp. $x_{20}$)
(resp. $x_{31}$) act on $V$ as  $K$ (resp. $K^*$)
 (resp. $K^{-1}$)
 (resp. $K^{*-1}$).
With reference to Corollary 
\ref{cor:four}
and 
Definition \ref{def:spltop}
we let $x_{01}$
(resp. $x_{12}$)
(resp. $x_{23}$)
(resp. $x_{30}$) act on $V$
as the split operator for $K,K^*$ 
(resp. $K^*,K^{-1}$)
(resp. $K^{-1},K^{*-1}$)
(resp. $K^{*-1},K$).
We now show that the above actions induce a $\boxtimes_q$-module
structure on $V$. To do this we show that they
 satisfy the defining relations
for $\boxtimes_q$ given in
Definition \ref{def:qtet}.
The relations in 
Definition \ref{def:qtet}(i) hold by the construction,
so consider the relations
in 
Definition \ref{def:qtet}(ii). The three kinds  
of relations involved are treated in the following
three claims.
\\

\noindent 
{\it
Line (\ref{eq:t2}) holds
for all $r,s,t \in \Z_4$ such that $(s-r,t-s)=(1,2)$.}
\hfil\break
\noindent 
To prove this claim, by
$\Z_4$ symmetry we may assume that
$r=0$, $s=1$, $t=3$.
Let $\Delta$ denote the
left hand side of
(\ref{eq:t2})
minus the right
hand side of 
(\ref{eq:t2}), for $r=0$, $s=1$, $t=3$.
We show $\Delta=0$. For $0 \leq i \leq d$ we combine
Lemma
\ref{lem:qweyl} (with $A=x_{01}$, $B=x_{13}$, $\theta=q^{d-2i}$)
and Lemma
\ref{lem:kks}(ii)
to find
$\Delta U_i=0$. Now $\Delta =0$ in view of Lemma
\ref{lem:udec} and the claim is proved.
\\

\noindent {\it 
Line (\ref{eq:t2}) holds
for all $r,s,t \in \Z_4$ such that $(s-r,t-s)=(2,1)$.}
\hfil\break
\noindent 
To prove this claim, by 
 $\Z_4$ symmetry we may assume that
 $r=2$, $s=0$, $t=1$.
Let $\Delta$ denote the
left hand side of
(\ref{eq:t2})
minus the right
hand side of 
(\ref{eq:t2}), for $r=2$, $s=0$, $t=1$.
We show $\Delta=0$. For $0 \leq i \leq d$ we combine
Lemma
\ref{lem:qweyl2} (with $A=x_{20}$, $B=x_{01}$, $\theta=q^{d-2i}$)
and Lemma 
\ref{lem:kks}(i)
to find
$\Delta U_i=0$. Now $\Delta =0$ in view of Lemma
\ref{lem:udec} and the claim is proved.
\\

\noindent {\it Line
 (\ref{eq:t2}) 
holds for 
all $r,s,t \in \Z_4$ such that $(s-r,t-s)=(1,1)$.}
\hfil\break
\noindent 
To prove the claim, by
$\Z_4$ symmetry it suffices to
show that 
(\ref{eq:t2}) holds
for $r=0$, $s=1$, $t=2$.
Combining 
(\ref{eq:t2}) (with 
$r=3$, $s=1$, $t=2$)
and Lemma
\ref{lem:qweyl}
we find
\begin{eqnarray*}
(x_{12}-q^{d-2i}I)V^*_i \subseteq V^*_{i-1} \qquad \qquad (0 \leq i\leq d).
\end{eqnarray*}
By this and Lemma
\ref{lem:uv}(ii)
we have
\begin{eqnarray}
(x_{12}-q^{2i-d}I)U_i \subseteq U_{i+1}+ \cdots + U_d
\qquad \quad (0 \leq i \leq d).
\label{eq:step1}
\end{eqnarray}
Combining 
(\ref{eq:t2}) (with 
$r=1$, $s=2$, $t=0$)
and Lemma \ref{lem:qweyl2}
we find
\begin{eqnarray*}
(x_{12}-q^{d-2i}I)V_i \subseteq V_{i+1} \qquad \qquad (0 \leq i\leq d).
\end{eqnarray*}
By this and Lemma
\ref{lem:uv}(i)
we have
\begin{eqnarray}
(x_{12}-q^{2i-d}I)U_i \subseteq U_0+ \cdots + U_{i+1}
\qquad \quad (0 \leq i \leq d).
\label{eq:step2}
\end{eqnarray}
Combining 
(\ref{eq:step1}), (\ref{eq:step2})
we find
\begin{eqnarray}
(x_{12}-q^{2i-d}I)U_i \subseteq U_{i+1}
\qquad \quad (0 \leq i \leq d).
\label{eq:step3}
\end{eqnarray}
Let $\Delta$ denote the left hand side of 
(\ref{eq:t2}) minus the right hand side of
(\ref{eq:t2}),
for $r=0$, $s=1$, $t=2$.
For $0 \leq i \leq d$ we 
combine
Lemma
\ref{lem:qweyl}  (with $A=x_{01}$, $B=x_{12}$, $\theta=q^{d-2i}$)
and 
(\ref{eq:step3})
to find
$\Delta U_i=0$.
Now $\Delta=0$ in view of Lemma
\ref{lem:udec}.
We have now shown that 
(\ref{eq:t2}) holds for
$r=0$, $s=1$, $t=2$, and the claim follows.
\\

\noindent 
We have now shown that
the relations in Definition \ref{def:qtet}(ii) are satisfied.
We now consider the relations in
 Definition \ref{def:qtet}(iii).
Let 
 $r,s,t,u \in \Z_4$ be given such that $s-r=t-s=u-t=1$.
We show that 
(\ref{eq:qserre}) holds.
By $\Z_4$ symmetry we may assume that
$r=0$, $s=1$, $t=2$, $u=3$.
Combining 
(\ref{eq:t2}) (with $r=0$, $s=2$, $t=3$)
and Lemma 
\ref{lem:qweyl} we find
\begin{eqnarray*}
(x_{23}-q^{2i-d}I)V_i \subseteq V_{i+1}
\qquad (0 \leq i \leq d).
\end{eqnarray*}
By this and Lemma
\ref{lem:uv}(i) we have 
\begin{eqnarray}
x_{23}U_i \subseteq U_0+\cdots + U_{i+1}
\qquad \qquad (0 \leq i \leq d).
\label{eq:x23act}
\end{eqnarray}
By (\ref{eq:t2}) (with $r=2$, $s=3$, $t=1$)
and Lemma
\ref{lem:qweyl2}
we find  
\begin{eqnarray*}
(x_{23}-q^{d-2i}I)V^*_i \subseteq V^*_{i+1}
\qquad (0 \leq i \leq d).
\end{eqnarray*}
By this and Lemma
\ref{lem:uv}(ii) we have  
\begin{eqnarray}
x_{23}U_i \subseteq U_{i-1}+\cdots + U_d
\qquad \qquad (0 \leq i \leq d).
\label{eq:x23act2}
\end{eqnarray}
Combining 
(\ref{eq:x23act}),
(\ref{eq:x23act2})
we find
\begin{eqnarray}
x_{23}U_i \subseteq U_{i-1}+U_i + U_{i+1}
\qquad (0 \leq i \leq d).
\label{eq:x23}
\end{eqnarray}
Let $\Delta$ denote the left hand side
of 
(\ref{eq:qserre}) for 
$r=0$, $s=1$, $t=2$, $u=3$.
For $0 \leq i \leq d$ we combine Lemma
\ref{lem:key} 
(with $A=x_{01}$, $B=x_{23}$, $\theta=q^{d-2i}$)
and 
(\ref{eq:x23})
 to find
$\Delta U_i=0$. Now 
$\Delta=0$ in view of
Lemma
\ref{lem:udec}.
We have now shown that
(\ref{eq:qserre}) holds
for $r=0$, $s=1$, $t=2$, $u=3$ as desired.
We have now verified all the relations in
Definition \ref{def:qtet}(iii).
\\

\noindent We have shown that the
actions on $V$ of the generators $x_{rs}$
satisfy the defining relations
for $\boxtimes_q$ given in Definition \ref{def:qtet}.
Therefore these actions induce a
$\boxtimes_q$-module structure
on $V$.
By the construction $x_{02}, x_{13}$ act
on $V$ as $K,K^*$ respectively.
Apparently there exists a 
$\boxtimes_q$-module structure
on $V$ such that 
$x_{02}, x_{13}$ act
on $V$ as $K,K^*$ respectively.
We now show that this $\boxtimes_q$-module structure
is unique.
 Suppose we are given
any $\boxtimes_q$-module structure on $
V$ such that 
$x_{02}, x_{13}$
act as $K, K^*$ respectively.
This $\boxtimes_q$-module structure is irreducible by construction
and Definition \ref{def:inv}(iii). 
This $\boxtimes_q$-module structure is 
type $1$ and diameter $d$, since the action
of $x_{01}$ on $V$ has eigenvalues $\lbrace q^{d-2i} \;|\;0 \leq i \leq d
\rbrace $. For each generator $x_{rs}$ of $\boxtimes_q$ the
action on $V$ is determined by the decomposition $\lbrack r,s\rbrack$.
By Theorem
\ref{thm:flagdec} the decomposition $\lbrack r,s\rbrack$
is determined by the flags  $\lbrack r \rbrack $ and
$\lbrack s \rbrack$.
Therefore our $\boxtimes_q$-module structure on $V$
is determined by the flags $\lbrack n\rbrack, n \in \Z_4$.
By construction the flags $\lbrack 0\rbrack$ and $\lbrack 2 \rbrack$
are determined by the decomposition $\lbrack 0,2 \rbrack $
and hence by the action of $x_{02}$ on $V$.
Similarly 
 the flags $\lbrack 1\rbrack$ and $\lbrack 3 \rbrack$
are determined by the decomposition $\lbrack 1,3 \rbrack $
and hence by the action of $x_{13}$ on $V$.
Therefore the given $\boxtimes_q$-module structure on $V$
is determined by the action of $K$ and $K^*$ on $V$,
so this 
$\boxtimes_q$-module structure is unique.
We have now shown that 
there exists
a unique
 $\boxtimes_q$-module
structure on $V$ such that
$x_{02}, x_{13}$
act as $K, K^*$ respectively.
We mentioned earlier that this
 $\boxtimes_q$-module
structure is irreducible and has type $1$.
\hfill $\Box $ \\

\section{Directions for further research}

\noindent In this section we give some
directions for further research. We start with a definition.

\begin{definition}
\label{def:cald}
\rm
Let $V$ denote a vector space over $\K$ with
finite positive dimension.
Let $\lbrace V_i\rbrace_{i=0}^d$ 
denote a decomposition of $V$.
For $0 \leq i \leq d$ let $F_i :V\to V$
denote the linear transformation that
satisfies
\begin{eqnarray*}
(F_i-I)V_i &=& 0, \\           
F_iV_j&=&0 \quad \mbox{if} \quad i\not=j
\qquad \qquad (0 \leq j\leq d).
\end{eqnarray*}
We observe that $F_i$ is the projection from $V$ onto
$V_i$. We note that $I=\sum_{i=0}^d F_i$ and
$F_iF_j=\delta_{ij}F_i$ for
$0 \leq i,j\leq d$. Therefore
the sequence $F_0,\ldots, F_d$ is a basis 
for a commutative subalgebra $\mathcal D$ of $\mbox{End}(V)$.
\end{definition}

\begin{problem}
\label{prob:genlp}
\rm
Let $V$ denote a vector space over $\K$ with
finite positive dimension.
Let $\lbrace V_i\rbrace_{i=0}^d$ 
and $\lbrace V^*_i\rbrace_{i=0}^{\delta}$ 
denote decompositions of $V$.
Let $\cal D$ and ${\cal D}^*$ 
denote the corresponding commutative algebras
from Definition
\ref{def:cald}. Investigate the case in
which (i)--(v) hold below.
\begin{enumerate}
\item[{\rm (i)}]  
$\cal D$ has a generator $A_{+}$ such that
\begin{eqnarray*}
A_{+} V^*_i \subseteq V^*_0+\cdots + V^*_{i+1}
\qquad (0 \leq i \leq \delta).
\end{eqnarray*}
\item[{\rm (ii)}]  
$\cal D$ has a generator $A_{-}$ such that
\begin{eqnarray*}
A_{-} V^*_i \subseteq V^*_{i-1}+\cdots + V^*_{\delta}
\qquad (0 \leq i \leq \delta).
\end{eqnarray*}
\item[{\rm (iii)}]  
${\cal D}^*$ has a generator $A^*_{+}$ such that
\begin{eqnarray*}
A^*_{+} V_i \subseteq V_0+\cdots + V_{i+1}
\qquad (0 \leq i \leq d).
\end{eqnarray*}
\item[{\rm (iv)}]  
${\cal D}^*$ has a generator $A^*_{-}$ such that
\begin{eqnarray*}
A^*_{-} V_i \subseteq V_{i-1}+\cdots + V_d
\qquad (0 \leq i \leq d).
\end{eqnarray*}
\item[{\rm (v)}]  
There does not exist a subspace $W\subseteq V$
such that ${\cal D}W\subseteq W$
and ${\cal D}^*W\subseteq W$,
other than $W=0$ and $W=V$.
\end{enumerate}
\end{problem}

\begin{note}
\rm 
Let $A,A^*$ denote a tridiagonal
pair on $V$, as in 
\cite[Definition 1.1]{TD00}. Then 
the conditions (i)--(v) of Problem
\ref{prob:genlp} are satisfied with
\begin{eqnarray*}
A_+=A, \qquad
A_-=A, \qquad
A^*_+=A^*, \qquad
A^*_-=A^*.
\end{eqnarray*}
\end{note}

\begin{note}
\rm 
Let $K,K^*$ denote a $q$-inverting pair on $V$.
Then 
the conditions (i)--(v) of Problem
\ref{prob:genlp} are satisfied with
\begin{eqnarray*}
A_+=K^{-1}, \qquad
A_-=K, \qquad
A^*_+=K^*, \qquad
A^*_-=K^{*-1}.
\end{eqnarray*}
\end{note}


\begin{problem} \rm
Referring to Problem
\ref{prob:genlp}, assume conditions (i)--(v) hold.
Show that the decompositions 
$\lbrace V_i\rbrace_{i=0}^d$ and
$\lbrace V^*_i\rbrace_{i=0}^\delta$ 
have the same shape and in particular that $d=\delta$.
Denoting this common shape by $\lbrace \rho_i\rbrace_{i=0}^d$,
show that $\rho_i = \rho_{d-i}$ for $0 \leq i \leq d$.
\end{problem}

\begin{problem} \rm
Referring to Problem
\ref{prob:genlp}, assume conditions (i)--(v) hold.
Consider the four decompositions of $V$ 
consisting of $\lbrace V_i\rbrace_{i=0}^d$ and
$\lbrace V^*_i\rbrace_{i=0}^\delta$, together with
their inversions.
Show that the flags on $V$ induced by these
four decompositions are mutually opposite.
\end{problem}

\noindent We motivate our last problem with
 a comment. By Theorem
\ref{thm:tdhard} and
(\ref{eq:qserre}) 
every $q$-tridiagonal
pair satisfies the cubic $q$-Serre relations.
More generally, every tridiagonal pair 
\cite[Definition 1.1]{TD00}
satisfies a 
pair of equations called the {\it tridiagonal relations}
\cite[Theorem 3.7]{qSerre}.

\begin{problem}
\rm Find some polynomial relations 
satisfied by every $q$-inverting pair.

\end{problem}

\noindent Tatsuro Ito \hfil\break
\noindent Department of Computational Science \hfil\break
\noindent Faculty of Science \hfil\break
\noindent Kanazawa University \hfil\break
\noindent Kakuma-machi \hfil\break
\noindent Kanazawa 920-1192, Japan \hfil\break
\noindent email:  {\tt ito@kappa.s.kanazawa-u.ac.jp}

\bigskip

\noindent Paul Terwilliger \hfil\break
\noindent Department of Mathematics \hfil\break
\noindent University of Wisconsin \hfil\break
\noindent Van Vleck Hall \hfil\break
\noindent 480 Lincoln Drive \hfil\break
\noindent Madison, WI 53706-1388 USA \hfil\break
\noindent email: {\tt terwilli@math.wisc.edu }\hfil\break

\end{document}